\documentclass[preprint,12pt]{article}
\usepackage{graphicx}
\usepackage[round,authoryear]{natbib}
\usepackage{amssymb}
\usepackage{amsmath}
\newtheorem{thm}{Theorem}
\newtheorem{lem}{Lemma}
\newtheorem{cor}{Corollary}
\newtheorem{prop}{Proposition}
\newcommand{\finish}{\hfill$\Box$\vspace{0.2cm}}

\newcommand{\prf}{\noindent{\bf Proof:\ }}

\newcommand{\E}{{\rm I \!E}}
\newcommand{\p}{{\rm I \!P}}

\providecommand{\keywords}[1]{\textbf{Keywords} #1}
\begin{document}

\title{Path Decomposition of Spectrally Negative L\'{e}vy Processes \thanks{This work is supported by Middle East Technical University BAP-08-11-2016-069}}



\author{Ceren Vardar-Acar \and Mine \c{C}a\u{g}lar}



\maketitle

\begin{abstract}  Path decomposition is performed to analyze the pre-supremum, post-supremum, post-infimum and the intermediate processes of a spectrally negative L\'evy process taken up to an independent exponential time $T$ as motivated by the aim of finding the joint distribution of the maximum loss and maximum gain. In addition, the joint distribution of the supremum and the infimum before an exponential time is displayed. As an application of path decomposition, the distributions of supremum of the post-infimum process and the maximum loss of the post-supremum process are obtained.
\end{abstract}

\section{Introduction}
\label{intro}

A spectrally negative L\'evy process $X$ is that with no positive jumps and has found place naturally in many applications such as risk theory, mathematical finance and queuing theory, e.g. \citep{Avra,baur,Kypr,MijoPist}. It is an advantageous model due to its tractability through the so-called scale functions.  In this paper, path decomposition of a spectrally negative L\'evy process is performed to analyze the pre-supremum, post-supremum, post-infimum and the intermediate processes as motivated by questions which can be expressed through trajectories and extreme values.

Let $X$ be defined on the filtered probability space $(\Omega, \mathcal {H},(\mathcal{ F}_t)_{t\ge 0}, \p)$ satisfying the usual hypotheses.
By definition, the L\'evy measure $\Pi$ of a spectrally negative L\'{e}vy process  is concentrated on $(-\infty, 0)$. Let the Laplace exponent   of $X$ be given by
\[
\psi (\lambda ) =  \mu \lambda  + \frac{{{\sigma ^2}}}{2}{\lambda ^2} + \int\limits_{( - \infty ,0)}^{} {({e^{\lambda x}} - 1 - \lambda x\;1_{\{x > -1\}})\, \Pi (dx)}
\]
where $\mu \in \mathbb{R}$, and $\sigma>0$, and let $X_{0}=0$. Then, Brownian motion with drift $\mu$ can be recovered as a special case when $\Pi\equiv 0$. Using a Brownian motion $B$ and a Poisson random measure $N$, for
$\mu \in \mathbb{R}$ and $\sigma >0$, we can write
\[
dX_t=\mu\, dt+\sigma\, dB_t + \int_{(-1,0)} y\, \widetilde{N}(dy,dt) + \int_{(-\infty,-1]}y\, N(dy,dt)
\]
where $\widetilde{N}=N-\Pi(dy)\, dt$. We assume $\int_ {-\infty}^0 (1 \wedge y^2 ) \, \Pi(dy)< \infty$.
The infinitesimal generator  of $X$  is given by
\[
 \mu f^\prime(x)+\frac{1}{2}\sigma^2  f^{\prime\prime}(x)+ \int_{-\infty}^0 [f(x+y)-f(x) -f^\prime(x)y\,1_{\{y>-1\}}]\, \Pi(dy)
\]
for  $f\in C_b^2$. We denote the first passage times above and below $x$ respectively by
\begin{equation*} \tau^{+}_x = \inf\{t\geq 0 : X_t>x\} ~~~~~~\tau^{-}_x = \inf\{t\geq 0 : X_t<x\} \;.\end{equation*}

We consider a spectrally negative L\'evy process $X$ taken up to an exponential time $T$ independent of $X$ with parameter $\gamma.$  Let the running supremum and infimum be denoted by $S$ and $I$, respectively, and let $S_{T}:=\sup \{ X_s: ~0  \le s \leq T\}$ and $I_{T}:=\inf \{ X_s: ~0  \le s \leq T\}$ before time $T$. The fluctuation identities that we rely on involve a class of functions known as scale functions. The $\gamma-$scale function of $X$ satisfies
\[\int\limits_0^\infty  {{e^{ - \lambda x}}} W^{(\gamma)}(x)dx = \frac{1}{{\psi (\lambda )}-\gamma} ~~~~\mbox{for}~~~~\lambda>\Phi(\gamma)  \]
where $\Phi$  denotes the right inverse of $\psi$ and the second scale function is given by
\begin{equation}  \label{Z}
{Z^{(\gamma)}}(x) = 1 + \gamma\int_0^x {{W^{(\gamma)}}(y)dy}
\end{equation}
see e.g. \citep[Thm. 8.1]{Kypr}.

  We perform path decomposition through the extremes and find the distributions of pre-supremum, post-supremum and post-infimum processes in Theorem \ref{thm3}. As a corollary, the cumulative distribution function of the supremum of the post-infimum process given $I_{t}=a$ turns out to be
$$\frac{\Phi(\gamma)({Z^{(\gamma)}}(b-a)-1)}{\gamma W^{(\gamma)}(b-a)}$$
 at $b\geq a$.
Then, Proposition \ref{prop2} characterizes the post-infimum process up to a given level. Its infinitesimal generator is given by
 \begin{eqnarray*}
{\mathcal L}F(x) & =&  \left[\mu +\sigma \frac{W^{ (\gamma) \prime}(x)-\Phi(\gamma)W^{ (\gamma)}(x)}
{W^{(\gamma)}(x)}\right] F'(x)+\frac{\sigma^2}{2} F''(x) \\
& & + \int_{(-1,0)}y \frac{e^{-\Phi(\gamma)y}W^{ (\gamma)}(x+y) - W^{ (\gamma)}(x)} {W^{(\gamma)}(x) }      \,\Pi(dy) \\
&& + \int_{-\infty}^0 [F(x+y)-F(x)-F'(x)\, y\, 1_{\{y>-1\}} ]     \\
&& \hspace{1cm} .\,\left[1+ \frac{e^{-\Phi(\gamma)y}W^{ (\gamma)}(x+y) - W^{ (\gamma)}(x)} {W^{(\gamma)}(x) } \right] \,\Pi(dy)
\end{eqnarray*}
for $F \in C^2_b$.
As in \citep[Sec.3]{tanre}, we apply the theory of expansion of filtrations to have the last exit from the infimum a stopping time and analyze the intermediate process between the infimum and the supremum.
On the other hand, a path decomposition through the infimum of  a L\'{e}vy process conditioned to stay positive has been studied by \citet{chaum}.
There, the focus is on a Levy processes conditioned to stay positive, whereas we focus on the spectrally negative Levy process killed at an exponential time $T$. Although related, the two processes and hence the results cannot be obtained from each other.

 In Theorem \ref{thm2}, we make use of Proposition \ref{prop2} under the condition that the infimum of the L\'evy process occurs before its supremum and at given levels to identify the law of the intermediate process between the infimum and the supremum. Under the same conditions, we also provide the distribution of the post-supremum process as a part of path decomposition. Such a path decomposition is motivated by an ultimate aim of finding the joint distribution of the maximum loss and the maximum gain defined by
\[ M_t^ - : = \mathop {\sup }\limits_{0 \le u \le v \le t} (X_u^{} - X_v^{})~~~~~~M_{t}^ +  := \mathop {\sup }\limits_{0 \le u \le v \le t} (X_v^{} - X_u^{})  \]
which are also known as maximum drawdown and maximum drawup. This approach has been originally followed by Salminen and Vallois \citep{salm} to obtain the same joint distribution for standard Brownian motion. In our case, the characterization of pre-infimum process remains unsolved. Nevertheless, the cumulative distribution function of the maximum loss of the post-supremum process at $0<d<b-a$ under the condition that infimum occurs before the supremum and $I_T=a$, $S_T=b$ is found as
\begin{eqnarray*}
  1- \left[1-Z^{(\gamma)}(b-a-d)+(Z^{(\gamma)}(b-a)-1)\frac{W^{ (\gamma)}(b-a-d)}{W^{(\gamma)}(b-a)}\right] \\
\cdot \,\frac{Z^{(\gamma)'}(d)-Z^{(\gamma)}(d)\frac{W^{ (\gamma)'}(d)}{W^{(\gamma)}(d)}}{-Z^{(\gamma)\prime}(b-a)+(Z^{(\gamma)}(b-a)-1)\frac{W^{(\gamma)'}(b-a)}{W^{(\gamma)}(b-a)}}
\end{eqnarray*}
as given in Corollary \ref{cor2}. We also provide the joint distribution of the infimum and the supremum of $X$ as $$\p_{0}(a<I_{T},S_{T}<b)=1-Z^{(\gamma)}(-a)
+[Z^{(\gamma)}(b-a)-1]\frac{W^{(\gamma)}(-a)}{W^{(\gamma)}(b-a)}$$
in terms of the scale functions. It can be combined with path decomposition to remove the conditions in our main results.

The paper is organized as follows.   The path of a spectrally negative L\'{e}vy process  is decomposed through its infimum and supremum separately before an independent exponential time in Section 2. The post-infimum process up to a given level is characterized in Section 3. Finally, Section 4 provides a more detailed path decomposition including the intermediate process between the infimum and the supremum, and displays the distribution of the maximum loss of the post-supremum process as an application.

\section{Decomposition through infimum or supremum}

In this section, we perform a path decomposition of the L\'evy process to characterize the distributions of pre-supremum, post-supremum, and post-infimum processes before an independent exponential time $T$ with parameter $\gamma>0$.
Let
\[
H_S:=\sup\{t < T: X_{t}=S_t\} \quad \mbox{and} \quad H_I:=\sup\{t < T: X_t=I_t\} \; .
\]
Then, we have the following theorem for pre-$H_S$ process $\{X_u:0\le u\le H_{S}\}$, post-$H_S$ process $\{X_{H_S+u}:0\le u\le T-H_S\}$, and post-$H_I$ process $\{X_{H_I+u}:0\le \emph{}u\le T-H_I\}$.

\begin{thm} \label{thm3} \begin{enumerate}
\item Conditionally on $I_T,$ the pre-$H_I$  process and the post-$H_I$ process are independent.
  Conditioned on $I_T$, the value of the infimum, the law of the post-$H_I$ process is given by $h$-transform of the law $\p^\uparrow$ of the original spectrally negative L\'evy process conditioned to stay positive and killed at time $T$ with
\[
 {h}(z)= \frac{\gamma \, W^{(\gamma)}(z)}{\Phi(\gamma)W(z)}-\frac{\gamma}{W(z)} \int_0^z W^{(\gamma)}(r)\, dr \; ,
\]
that is, the transition semigroup of the post-$H_I$ process when $I_T=a$ is given by
\[
P_t(x,dy)=\frac{h(y-a)}{h(x-a)}\, \p^\uparrow _{x-a}  \left\{ X_t\in dy-a, t<T \right\}
\]
for $x,y>a$, and its entrance law is obtained by letting $x\to a+$. \\
\item  Conditionally on $S_T$, the pre-$H_S$  process and the post-$H_S$ process are independent. Conditioned on $S_T$, the value of the supremum, the pre-$H_S$ process is a spectrally negative L\'evy process with Laplace exponent
\[
\bar{\psi}(\lambda)=\psi(\lambda + \Phi(\gamma))-\gamma
\]
for $\lambda\ge -\Phi(\gamma)$, killed at the first passage time above $b$ when $S_T=b$. On the other hand, when $S_T=b$, the transition semigroup of the post-$H_S$ process is given by
\[
P_t(x,dy)=\frac{h(b-y)}{h(b-x)}\, \p_{b-x}  \left\{ \hat{X}_t\in b-dy, t<T\wedge \hat{\tau}_0^- \right\}
\]
for $x,y<b$, where
\[
h(z)=1-e^{-\Phi(\gamma)z}\, ,
\]
$\hat{X}$ is a spectrally positive L\'{e}vy process distributed like $-X$, and the entrance law is obtained as $x\to b-$.
\end{enumerate}
\end{thm}

\prf
\textbf{i.} The independence of pre-$H_I$ and post-$H_I$ process conditionally on $I_T$ follows from \citep[Lem.VI.6.ii]{bert} by applying it to the process $-X$ noting that $0$ is regular for $(0,\infty)$ since $X$ is assumed to be of unbounded variation. Also note that $X_{H_I}=I_T$ almost surely since we have assumed $\sigma>0$ and the process creeps downwards by \citet[Ex.7.6.iv]{Kypr}. The post-$H_I$ process is a Markov process with stationary transitions, which depend on the value of the infimum, $I_T$ \citep{mill}.
The post-$H_I$ process evolves as a spectrally negative L\'evy process which is conditioned to stay above $I_T$ and killed at an exponential time $T$. We will identify its transition function $P_t$ by studying
\begin{equation}
 P_t(x,dy) := \p_x\{X_t\in dy, t<T\,|\, T< \tau_0^-\}  \label{post}
\end{equation}
which is equal to
\begin{eqnarray*}
& &  \displaystyle{\frac{  \p_x\{X_t\in dy, t<T \wedge \tau_0^-, T< \tau_0^-\} }{\p _x\{T< \tau_0^-\} }} \nonumber \\
& & =\displaystyle{ \frac{  \p_x\{ T< \tau_0^- | X_t = y, t<T \wedge \tau_0^- \}\p _x  \left\{   X_t\in dy, t<T \wedge \tau_0^- \right\}}{\p _x\{T< \tau_0^-\} } } \nonumber \\
& & = \displaystyle{\frac{  \p_y\{ T\circ \theta_t< \tau_0^-\circ \theta_t \}\p _x  \left\{   X_t\in dy, t<T \wedge \tau_0^- \right\}}{\p _x\{T< \tau_0^-\} } }\nonumber \\
& & = \displaystyle{ \frac{  \p_{y} \{ T< \tau_0^-\} }{\p _x\{T< \tau_0^-\}} \p _x  \left\{   X_t\in dy, t<T \wedge \tau_0^- \right\} } \nonumber
\end{eqnarray*}
by Markov property and the fact that $T$ is memoryless, where $\theta$ is the shift operator \citep[Sec.7.3]{cinl}. Let
\begin{equation}\tilde{h}(x) := \p _x\{T< \tau_0^-\} \; .\nonumber\end{equation}
The function $\tilde{h}$ can be evaluated as
\begin{eqnarray}
\tilde{h}(x) &=& \E_x[\p_x\{T<\tau_0^-|\, \tau_0^-\}] = \E_x[1-e^{-\gamma \tau_0^-}1_{\{\tau_0^-<\infty\}}] \nonumber \\
&=& 1-\E_x[e^{-\gamma \tau_0^-}1_{\{\tau_0^-<\infty\}}] \nonumber \\
&=& 1- Z^{(\gamma)}(x)+\frac{\gamma}{\Phi(\gamma)}W^{(\gamma)}(x)\nonumber \\
& = & \frac{\gamma W^{(\gamma)}(x)}{\Phi(\gamma)}-\gamma \int_0^x W^{(\gamma)}(z)\, dz \label{th}
\end{eqnarray}
by \citep[Thm.8.1]{Kypr} and (\ref{Z}).
Note that the post-$H_I$ process can now be viewed as a transform of the law of the L\'evy process conditioned to stay positive \citep[pg.198]{bert} by rearranging (\ref{post}) as
\begin{eqnarray*}
\lefteqn{\p_x\{X_t\in dy, t<T\,|\, T< \tau_0^-\}} \\
&  \displaystyle{=\frac{\tilde{h}(y)}{\tilde{h}(x) } \p _x  \left\{   X_t\in dy, t< \tau_0^-, t<T \right\} \quad \quad \quad \quad} \\
&  \displaystyle{ \quad =\frac{\tilde{h}(y)/W(y)}{\tilde{h}(x)/W(x) } e^{-\gamma t} \, \frac{W(y)}{W(x)} \, \p _x  \left\{   X_t\in dy, t< \tau_0^-\right\}   } \\
& \displaystyle{\quad \quad \quad \quad=:\frac{{h}(y)}{{h}(x) } e^{-\gamma t} \, \p^\uparrow _x  \left\{   X_t\in dy\right\} = \frac{{h}(y)}{{h}(x) } \, \p^\uparrow _x  \left\{   X_t\in dy, t<T \right\} }
\end{eqnarray*}
where $\p^\uparrow _x $ denotes the law of the L\'evy process started at $x$ and conditioned to stay positive, and $W(x)=W^{(0)}(x)$. Explicitly, $h=\tilde{h}/W$ is given by
\[
 {h}(x)= \frac{\gamma \, W^{(\gamma)}(x)}{\Phi(\gamma)W(x)}-\frac{\gamma}{W(x)} \int_0^x W^{(\gamma)}(z)\, dz \; .
\]
As a result, the entrance law of the post-$H_I$ process is ${h}$-transform of the law of the L\'evy process conditioned to stay positive and killed at an exponential time \citep[Prop.VII.14]{bert}. \\

\textbf{ii.} Note that $X_{H_{S}-}=S_T$ as $X$ is spectrally negative. Since we have assumed that $X$ is of unbounded variation, which in turn is equivalent to regularity of $0$  for $(-\infty,0)$, Lem.VI.6.ii of \citep{bert} for Markov processes yields the independence conditionally on $S_T$.
Considering the excursions of $X$ from its supremum $S$ as in the proof of this lemma, in particular for a spectrally negative L\'evy process, and with similar arguments as in the proof of Thm.3.2 of \citep{salm}, we derive the distribution of pre-$H_S$ process. Since $S$ serves as a local time for the reflected process $S-X$, we can express its excursions using $\tau$, which is a subordinator \cite[Thm.VII.1]{bert}. The points of discontinuity of $\tau=\{\tau_x:x\geq 0\}$ indicate the start of excursions $\varepsilon_x$. Explicitly, we have
\[
 \varepsilon_x = \left\{ {X_{{\tau_{{x }}}}} - {X_{{\tau_{x-}} + u}}{\kern 1pt} ,{\kern 1pt} \,0 < u \le \tau_{x} - \tau_{x-}  \right\}.
 \]
for $x>0$ such that $\tau_{x} - \tau_{x-}>0$ \cite{KuztKypr}.
If $\tau_{x} - \tau_{x-}=0$, take  $\varepsilon_x=\partial$ for some cemetery state $\partial$. Let ${\cal E}$ denote the space of real valued right continuous paths with left limits killed at the first hitting time of $(-\infty,0]$, endowed with the $\sigma$-algebra generated by coordinates. Then, $\{(x,\varepsilon_x):x> 0, \varepsilon_x \neq \partial\}$ constitute the atoms of a Poisson random measure $M$ on $(0,\infty)\times {\cal E}$ with  mean measure $dx \, n(d\varepsilon)$ for some  measure $n$ on ${\cal E}$; see \cite{bert} for a construction of $n$. Let the lifetime of a generic excursion be defined by ${\zeta}= \inf \left\{ u > 0:\varepsilon (u) \le 0 \right\}$. For bounded measurable functionals $F_1$ and $F_2$, we consider the following expectation to identify the distribution of the pre-supremum process
\begin{eqnarray} \label{product}
\lefteqn{\displaystyle{\E\left[ F_1(X_u:u\le H_{S})F_2(S_T-X_{H_S+u}:u\le T-H_S)\right] } } \\
&  \displaystyle{\!\!\!\!\! = \E\left[\int_0^\infty\!\! \int_{\mathcal{ E}} \!\! M(dx, d\varepsilon_x) \, F_1(X_u:u\le \tau_x)\,1_{\{\tau_x<T\}} F_2 (\varepsilon_x(u): u\le T-\tau_x)\, 1_{\{T-\tau_x< \zeta\}}(\varepsilon_x) \right]
}  \nonumber\\
& \displaystyle{\!\!\!\!\! = \E\left[\int_0^\infty\!\! \int_{\mathcal{ E}} \!\! dx\,  n(d\varepsilon) \, F_1(X_u:u\le \tau_x)\,1_{\{\tau_x<T\}} F_2 (\varepsilon(u): u\le T)\, 1_{\{T< \zeta\}}(\varepsilon) \right]
}  \nonumber\\
& \displaystyle{ \!\!\!\!\! = \E\left[\int_0^\infty \!\! dx \, F_1(X_u:u\le \tau_x)\,1_{\{\tau_x<T\}} \int_{\mathcal{ E}}\! n(d\varepsilon) F_2 (\varepsilon(u): u\le T)\, 1_{\{T< \zeta\}}(\varepsilon) \right]} \nonumber\\
& \displaystyle{ \!\!\!\!\! = \E\left[\int_0^\infty \!\! dx \, F_1(X_u:u\le \tau_x)\,e^{-\gamma \tau_x}\right] \E \left[\int_{\mathcal{ E}}\! n(d\varepsilon) F_2 (\varepsilon(u): u\le T)\, 1_{\{T< \zeta\}}(\varepsilon) \right]} \nonumber
\end{eqnarray}
where the steps follow from Markov property and
  the fact that $T$ is an independent exponential random variable  with parameter $\gamma>0$.

For  the law of the pre-$H_S$ process, by changing $x$ to $b$,\emph{} we observe that
\begin{eqnarray*}
\lefteqn{\displaystyle{\E [ F_1(X_u:u\le \tau_b)\,e^{-\gamma \tau_b}]}} \\
& \hspace{-4.7cm}\displaystyle{ = \int \p (d\omega) \, F_1(X_u(\omega):u\le \tau_b(\omega))\, e^{-\gamma \tau_b(\omega)} }\\
&  \displaystyle{ = \int \p^{\Phi(\gamma)}(d\omega) \, e^{-\Phi(\gamma)X_{\tau_b(\omega)}(\omega)+\psi(\Phi(\gamma))\tau_b(\omega)} F_1(X_u(\omega):u\le \tau_b(\omega))\, e^{-\gamma \tau_b(\omega)} }\\
 & \hspace{-6.2cm} \displaystyle{= \E^{\Phi(\gamma)} [F_1(X_u:u\le \tau_b)]\,e^{-\Phi(\gamma)b}}
\end{eqnarray*}
by change of measure with $c=\Phi(\gamma)$ as given in \citep[pg.232]{Kypr}, where $\p^c$ is the law of another spectrally negative L\'evy process with Laplace exponent $\bar{\psi}(\lambda)= \psi(\lambda+c)-\psi(c)$, and we put $X_{\tau_b} =b $ since the spectrally negative L\'{e}vy process $X$ creeps upwards almost surely. From (\ref{product}), we get
\begin{eqnarray} \label{joint}
\lefteqn{\displaystyle{\E\left[ F_1(X_u:u\le H_{S})F_2(S_T-X_{H_S+u}:u\le T-H_S)\right] } } \\
& \displaystyle{=\int_0^\infty \!\! db \, \E^{\Phi(\gamma)} [ F_1(X_u:u\le \tau_b)] \,e^{-\Phi(\gamma)b}  \int_{\mathcal{ E}}\! n(d\varepsilon)\, \E [F_2 (\varepsilon(u): u\le T)\, 1_{\{T< \zeta\}}(\varepsilon)]} \nonumber  \\
& \displaystyle{=\frac{1}{\Phi(\gamma)} \int_0^\infty \!\! db \, \E^{\Phi(\gamma)} [ F_1(X_u:u\le \tau_b)] \,f_{S_T}(b)  \int_{\mathcal{ E}}\! n(d\varepsilon)\, \E [F_2 (\varepsilon(u): u\le T)\, 1_{\{T< \zeta\}}(\varepsilon)]} \nonumber
\end{eqnarray}
where we wrote $f_{S_T}$ for the probability density function of $S_T$,
which has exponential distribution with parameter $\Phi(\gamma)$ \citep[Cor. VII.2]{bert}. The assertion about the conditional distribution of the pre-$H_S$ process given $S_T=b$ follows from (\ref{joint}).

The law of the post-$H_S$ process is characterized in view of  the dual process $\hat{X}:=-X$, which is spectrally positive. Similar to the analysis for post-$H_I$ process, we get
\begin{eqnarray}
& & \hspace{-1cm}\p_{b-x}\{\hat{X}_t\in b-dy, t<T\,|\, T<  \hat{\tau}_0^-\} \nonumber\\
& = & \frac{  \p_{b-y} \{T<\hat{\tau}_0^-\} }{\p _{b-x}\{T<  \hat{\tau}_0^-\}} \p_{b-x}  \{  \hat{X}_t\in b-dy, t<T \wedge  \hat{\tau}_0^- \} \nonumber
\end{eqnarray}
by Markov property and the fact that $T$ is memoryless.
This is the $h$-transform of the law of a spectrally positive L\'{e}vy process killed at the minimum of an exponential time and passage time below 0, with
\begin{eqnarray*}
h(x) &=& \p_x\{T< \hat{\tau}_0^-\}= \p_{-x} \{T<\tau_0^+\}=\p \{T<\tau_x^+\} \\
& =& 1-\E[e^{-\gamma \tau_x^+}1_{\{\tau_x^+<\infty\}}]\\
&=& 1-e^{-\Phi(\gamma)x}
\end{eqnarray*}
for $x>0$ \citep[pg.233]{Kypr}.
\finish

As a special case, the auxiliary function $\tilde{h}(x)$ as well as $h(x)$ corresponding to post-$H_S$ in the proof above are equal to $1-e^{-x\sqrt{2\gamma}}$ for standard Brownian motion due to symmetry and we recover the result in \citep[Thm.3.2.1]{salm}. Moreover, we obtain the following corollary, as a generalization of \citep[Eq.3.23]{salm}.
\begin{cor}
The distribution of the supremum of the post-infimum process, $S_{H_{I},T}$, when $I_T=a$ is given by
$$\p(S_{H_{I},T}\le b|I_T=a)=\frac{\Phi(\gamma)({Z^{(\gamma)}}(b-a)-1)}{\gamma W^{(\gamma)}(b-a)}$$
for $b>a$.
\end{cor}
\prf Note that the transition density of the post-infimum process can be written as
\[
P_t(x,dy)=\frac{\tilde{h}(b-a)}{\tilde{h}(x-a)}\p_{x-a}\left\{ X_t\in dy-a, t<T\wedge\tau_0^- \right\}
\]
for $x,y>a$ where $\tilde{h}$ is given in \eqref{th}. We have
\begin{eqnarray*} \p(S_{H_{I},T}>b|I_T=a)&=& \lim_{x \to a} \frac{\tilde{h}(b-a)}{\tilde{h}(x-a)}\, \p_{x-a}  \left\{ \tau_{b-a}^+<\tau_0^-, \tau_{b-a}^+<T \right\}\\
&=&\lim_{x \to a} \frac{\tilde{h}(b-a)}{\tilde{h}(x-a)}\frac{W^{(\gamma)}(x-a)}{W^{(\gamma)}(b-a)}\\
&=&\lim_{x \to a}\frac{\frac{\gamma W^{(\gamma)}(b-a)}{\Phi(\gamma)}-\gamma \int_0^{b-a} W^{(\gamma)}(z)\, dz }{\frac{\gamma W^{(\gamma)}(x-a)}{\Phi(\gamma)}-\gamma \int_0^{x-a} W^{(\gamma)}(z)\, dz }\, \frac{W^{(\gamma)}(x-a)}{W^{(\gamma)}(b-a)}\\
&=&1-\frac{\Phi(\gamma)}{W^{(\gamma)}(b-a)}\int_0^{b-a} W^{(\gamma)}(z)\, dz
\end{eqnarray*}
which is simplified by \eqref{Z} and we have used the fact that $W^{(\gamma)}(0)=0$ because $X$ is assumed to be of unbounded variation.
\finish
\section{Post-infimum process up to a given level}

In this section, we decompose the process according to its infimum until the first passage time above a given level which occurs before an independent, exponentially distributed time $T.$ The last exit time from the infimum is not a stopping time under the natural filtration generated by the spectrally negative L\'evy process, denoted by ${\cal F}$ below. On the other hand, after applying the theory of expansion of filtrations described in \citep[Sec.VI.3]{prot}, the last time the process exits from its infimum becomes a stopping time. As in \citep[Sec.3]{tanre}, we expand the filtration and characterize the distribution of the post-infimum process up to a given level before the exponential time.

For $b>0,$ we consider the trace of the original probability space on the event $D=\{ \tau_b^+ \le T\}$ \citep[Exer.I.1.15,I.3.12]{cinl}. On this probability space, define the random variable
\begin{equation}\label{rho}
 \rho = \sup \{ u : 0\le u\le \tau_b^+ \le T, X_u= I_{\tau^{+}_{b}}\}
 \end{equation}
as  the last exit time from the infimum before the first passage time above $b$, which occurs before an exponential time $T$. Here, $\rho$ is the end of the random set
\[
\Lambda= \{(u,\omega) \subset \mathbb{R}_+\times \Omega: 0\le u\le \tau_b^+(\omega) \le T(\omega), X_u(\omega)= I_{ \tau^{+}_{b}(\omega)} \} \; ,
\]
that is, $\rho= \sup \{ u  : (u, \omega) \in  \Lambda \}$, and we can apply the progressive expansions described  in \citep[Sec.VI.3]{prot}. The trace of the original filtration $\mathcal{F}$ is given by  $\mathcal{F}^{'}_{t}:=\mathcal{F}_{t}\cap \{\tau^{+}_{b}\leq T\}= \{A\cap D: A\in \mathcal{F}_t\} $.
Let $\mathcal{F}^{\rho}$ denote the smallest expanded filtration making $\rho$ a stopping time, satisfying the usual hypotheses. Let
\begin{equation}\label{yt}
Y_t = \p\{\rho\leq t | \mathcal{F}^{'}_{t}\} \quad \quad  t\geq 0\; .
\end{equation}
\begin{lem} \label{Y}  The process $Y$  is a submartingale with respect to $\mathcal{F}^{'}_{t}$ given by
 \[
 Y_{t}=\frac{e^{-\Phi(\gamma)X_{t\wedge\tau_b^+}}W^{(\gamma)}(X_{t\wedge \tau_b^+ }-I_{t\wedge \tau_b^+})}{e^{-\Phi(\gamma)b}W^{(\gamma)}(b-I_{t\wedge \tau_b^+})}
 \]
 for $t\geq 0$ with  Doob-Meyer decomposition  $Y=M+A$ where
 \begin{eqnarray*}
 \lefteqn{M_{t}
 = \int_0^{t\wedge \tau_b^+} \frac{e^{-\Phi(\gamma)X_{s}}[W^{ (\gamma) \prime}(X_s-I_{s})-\Phi(\gamma)W^{ (\gamma)}(X_{s}-I_{s})]}{e^{-\Phi(\gamma)b}W^{(\gamma)}(b -I_s)}\, \sigma d B_s  }  \label{mart} \\
 &
\displaystyle{ + \int_0^{t\wedge \tau_b^+} \int_{-\infty}^0  \frac{e^{-\Phi(\gamma)(X_{s^-}+y)}W^{ (\gamma)}(X_{s^-}+y-I_{s^-}) - e^{-\Phi(\gamma)X_{s^-}}W^{ (\gamma)}(X_{s^-}-I_{s^-})} {e^{-\Phi(\gamma)b}W^{(\gamma)}(b-I_{s^-})}  \,\widetilde{N}(dy,ds) }
 \end{eqnarray*}
is a martingale   and
 \begin{eqnarray*}
 {A}_t &=&\int_0^{t\wedge \tau_b^+ }  \frac{e^{-\Phi(\gamma)X_{s}}[W^{ (\gamma) \prime}(X_s-I_{s})-\Phi(\gamma)W^{ (\gamma)}(X_{s}-I_{s})]}{e^{-\Phi(\gamma)b}W^{(\gamma)}(b -I_s)}\,\mu ds\\
&+& \int_0^{t\wedge \tau_b^+ } \frac{\sigma^2}{2} \frac{e^{-\Phi(\gamma)X_{s}\left[\Phi^{2}(\gamma)W^{ (\gamma)}(X_s-I_s)-2\Phi(\gamma)W^{ (\gamma) \prime}(X_s-I_s)+W^{ (\gamma) \prime \prime}(X_s-I_s)\right]}}{e^{-\Phi(\gamma)b}W^{(\gamma)}(b -I_s)} ds \\
&+& \int_0^{t\wedge \tau_b^+ }\int_{-\infty}^0 \frac{e^{-\Phi(\gamma)(X_{s}+y)}W^{ (\gamma)}(X_{s}+y-I_{s}) - e^{-\Phi(\gamma)X_{s}}W^{ (\gamma)}(X_{s}-I_{s})}{e^{-\Phi(\gamma)b}W^{(\gamma)}(b -I_{s})}\,\Pi(dy) \, ds \\
&-& \int_0^{t\wedge \tau_b^+ }\int_{-\infty}^0 \frac{e^{-\Phi(\gamma)X_{s}}\left[W^{ (\gamma) \prime}(X_{s}-I_{s})-\Phi(\gamma)W^{ (\gamma)}(X_{s}-I_{s})\right]\, y \, 1_{\{y>-1\}}}{e^{-\Phi(\gamma)b}W^{(\gamma)}(b -I_{s})}\,\Pi(dy) \, ds \\
& &   + \int_0^{t\wedge \tau_b^+} \frac{e^{-\Phi(\gamma)X_{s}}[W^{ (\gamma) }(X_s-I_{s})W^{ (\gamma) \prime}(b-I_{s}) - W^{ (\gamma) \prime}(X_s-I_{s})W^{(\gamma)}(b -I_s)]}{e^{-\Phi(\gamma)b}[W^{(\gamma)}(b -I_s)]^2}\, d I_s^c \\
& & +\sum_{s\le t } \left[ \frac{e^{-\Phi(\gamma)X_{s}}W^{(\gamma)}(X_s-I_{s})}{e^{-\Phi(\gamma)b}W^{(\gamma)}(b-I_{s})} - \frac{e^{-\Phi(\gamma)X_{s^-}}W^{(\gamma)}(X_{s^-}-I_{s^-})}{e^{-\Phi(\gamma)b}W^{(\gamma)}(b-I_{s^-})} \right]
\end{eqnarray*}
is a nondecreasing predictable process.
\end{lem}

\prf
Let $t \geq 0$ be fixed.
When $\rho\le t$, after time $t$, the process passes above the level $b$ before it passes below the level $ I_{\tau^{+}_{b}}$, which is equal to $I_t$ in this case.
Therefore, we get
\[
\begin{split}
\p\{\rho\le  t| \mathcal{F}^{'}_{t}\}\, 1_{\{t\le  \tau^{+}_{b}\}}&=\frac{\p_{X_t}\{\tau^{+}_{b}< \tau^{-}_{I_{t}}, ~\tau^{+}_{b} \le T  \}}{\p_{X_t}\{\tau^{+}_{b} \le T \}} \, 1_{\{t\le  \tau^{+}_{b}\}}\nonumber\\
&=\frac{\E_{X_t}[e^{-\gamma \tau^{+}_{b}}1_{\{\tau^{+}_{b}<\tau^{-}_{I_{t}}\}}] }{\E_{X_t}[e^{-\gamma \tau^{+}_{b}}]} \, 1_{\{t\le  \tau^{+}_{b} \}}\nonumber\\
&=\frac{e^{-\Phi(\gamma)X_t}W^{(\gamma)}(X_t-I_{t})}{e^{-\Phi(\gamma)b}W^{(\gamma)}(b-I_{t})} \, 1_{\{t\le  \tau^{+}_{b} \}}
\end{split}
\]
 where $W^{(\gamma)}$ is the ${\gamma}$-scale function and the last equality follows from \citep[Thm.8.1]{Kypr}.
Therefore, we get
 \[
 Y_{t}=\frac{e^{-\Phi(\gamma)X_t}W^{(\gamma)}(X_{t}-I_{t})}{e^{-\Phi(\gamma)b}W^{(\gamma)}(b-I_{t})} \quad \quad t \leq \tau^{+}_{b} \; .
 \]
Note that for $t \geq \tau^{+}_{b}$ $Y_{t}=1$ by definitions \eqref{rho}, \eqref{yt}, and the expression for $Y_t$ follows. The process $Y_t$ is an ($\mathcal{F}^{'}_{t}$)-submartingale as it is immediate from definition \eqref{yt}.
By Ito's formula \citep[Thm.I.4.57]{jaco}, we have
\begin{eqnarray*}
d {Y}_t &=&\frac{e^{-\Phi(\gamma)X_{t^-}}[W^{ (\gamma) \prime}(X_{t^-}-I_{t^-})-\Phi(\gamma)W^{ (\gamma)}(X_{t^-}-I_{t^-})]}{e^{-\Phi(\gamma)b}W^{(\gamma)}(b -I_{t^-})} (\mu\, dt+\sigma d B_{t} ) \\
&&\hspace{-2cm}+ \frac{1}{2} \frac{e^{-\Phi(\gamma)X_{t^-}}[\Phi^{2}(\gamma)W^{ (\gamma)}(X_{t^-}-I_{t^-})-2\Phi(\gamma)W^{ (\gamma) \prime}(X_{t^-}-I_{t^-})+W^{ (\gamma) \prime \prime}(X_{t^-}-I_{t^-})]}{e^{-\Phi(\gamma)b}W^{(\gamma)}(b -I_{t^-})}\, \sigma^2 \, dt \\
&+& \int_{-\infty}^0 \frac{e^{-\Phi(\gamma)(X_{t^-}+y)}W^{ (\gamma)}(X_{t^-}+y-I_{t^-}) - e^{-\Phi(\gamma)X_{t^-}}W^{ (\gamma)}(X_{t^-}-I_{t^-})}{e^{-\Phi(\gamma)b}W^{(\gamma)}(b -I_{t^-})}\,\Pi(dy) \, dt \\
&-& \int_{-\infty}^0 \frac{e^{-\Phi(\gamma)X_{t^-}}\left[W^{ (\gamma) \prime}(X_{t^-}-I_{t^-})-\Phi(\gamma)W^{ (\gamma)}(X_{t^-}-I_{t^-})\right]\, y \, 1_{\{y>-1\}}}{e^{-\Phi(\gamma)b}W^{(\gamma)}(b -I_{t^-})}\,\Pi(dy) \, dt \\
& & + \int_{-\infty}^0 \left[ \frac{e^{-\Phi(\gamma)(X_{t^-}+y)}W^{ (\gamma)}(X_{t^-}+y-I_{t^-}) - e^{-\Phi(\gamma)X_{t^-}}W^{ (\gamma)}(X_{t^-}-I_{t^-})}{e^{-\Phi(\gamma)b}W^{(\gamma)}(b -I_{t^-})} \right]\,\widetilde{N}(dy,dt) \\
& &   + \frac{e^{-\Phi(\gamma)X_{t^-}}[W^{ (\gamma) }(X_{t^-}-I_{t^-})W^{ (\gamma) \prime}(b-I_{t^-}) - W^{ (\gamma) \prime}(X_{t^-}-I_{t^-})W^{(\gamma)}(b -I_{t^-})]}{e^{-\Phi(\gamma)b}[W^{(\gamma)}(b -I_{t^-})]^2}\, d I_t^c \\
& & +\frac{e^{-\Phi(\gamma)X_{t}}W^{(\gamma)}(X_t-I_{t})}{e^{-\Phi(\gamma)b}W^{(\gamma)}(b-I_{t})} - \frac{e^{-\Phi(\gamma)X_t}W^{(\gamma)}(X_{t^-}-I_{t^-})}{e^{-\Phi(\gamma)b}W^{(\gamma)}(b-I_{t^-})}
\end{eqnarray*}
where $I^c$ denotes the continuous part of the decreasing process $I$.

Now, $M_t$ is a martingale because both summands in its expression are martingales. First, the integrand of the Brownian integral is bounded as the process $X_t$ and $X_t-I_t$ is bounded when $t \in [0,\tau_b^+]$. Its quadratic variation is well defined because $W^{(\gamma)\prime}$ and $W^{(\gamma)\prime \prime}$ are both locally bounded at $0^{+}$ see e.g. \citep[pg.11]{avramvardar}. Then, the numerator  of the integrand in the Poisson integral can be bounded by a multiple of $y$ by mean value theorem as $W^{(\gamma)}$ is continuously differentiable and $X_t-I_t$ is bounded over $t \in [0,\tau_b^+]$. As a result, the second integral is a stochastic integral with respect to the compensated Poisson random measure under the conditions that the spectrally negative L\'{e}vy process $X$ exists. Its construction is similar to that of $X$, see e.g. \citep[Thm.2.10]{Kypr}. Therefore, the second term in $A_t$ is  well-defined. It is non-decreasing because $Y$ is a submartingale. Clearly, the other terms are well-defined as well. \finish

We can now characterize the post-$\rho$ process ${X_{t+\rho} - X_\rho}$  in the following proposition. Recall that $\rho = \sup \{ u : 0\le u\le \tau_b^+ \le T, X_u= I_{\tau^{+}_{b}}\}$.

\begin{prop}\label{prop2} Conditionally on $I_{\tau_b^+}=a$ and $\tau_b^+<T$,  the  process ${X_{\rho+t} - X_\rho}$ is a jump-diffusion process with generator
 \begin{eqnarray*}
{\mathcal L}F(x) & =&  \left[\mu +\sigma \frac{W^{ (\gamma) \prime}(x)-\Phi(\gamma)W^{ (\gamma)}(x)}
{W^{(\gamma)}(x)}\right] F'(x)+\frac{\sigma^2}{2} F''(x) \\
& & + \int_{(-1,0)}y \frac{e^{-\Phi(\gamma)y}W^{ (\gamma)}(x+y) - W^{ (\gamma)}(x)} {W^{(\gamma)}(x) }      \,\Pi(dy) \\
&& + \int_{-\infty}^0 [F(x+y)-F(x)-F'(x)\, y\, 1_{\{y>-1\}} ]     \\
&& \hspace{1cm} .\,\left[1+ \frac{e^{-\Phi(\gamma)y}W^{ (\gamma)}(x+y) - W^{ (\gamma)}(x)} {W^{(\gamma)}(x) } \right] \,\Pi(dy)
\end{eqnarray*}
for $F \in C^2_b$, starting at 0 and stopped at $b-a$.
\end{prop}

\prf  By \citep[Thm.VI.18]{prot}, if $C$ is a martingale for the original filtration then $C-\Upsilon$ is a martingale for the expanded filtration ${\cal F}^\rho$ where
\[
\Upsilon_t := -\int_0^{t\wedge \rho} \frac{1}{1-Y_s} d\langle C,M\rangle_s + 1_{\{t\ge \rho\}} \int_\rho ^t \frac{1}{Y_s} d\langle C,M\rangle_s
\]
and $Y$   and $M$ are as in Lemma \ref{Y}. We will consider the martingales that are relevant for $X$.

First consider the Brownian motion $B$. For $t< \tau_b^+$, define
\begin{eqnarray*}
\Lambda_t&=& -\int_0^{t\wedge \rho} \frac{1}{1-Y_s} d\langle B,M\rangle_s + 1_{\{t\ge \rho\}} \int_\rho ^t \frac{1}{Y_s} d\langle B,M\rangle_s \\
&=& -\int_0^{t\wedge \rho}    \frac{e^{-\Phi(\gamma)X_{s}}[W^{ (\gamma) \prime}(X_s-I_{s})-\Phi(\gamma)W^{ (\gamma)}(X_{s}-I_{s})]}{e^{-\Phi(\gamma)b}W^{(\gamma)}(b-I_{s})-e^{-\Phi(\gamma)X_{s}}W^{(\gamma)}
(X_{s}-I_{s})}\, \sigma \, ds   \\& &  + 1_{\{t\ge \rho\}} \int_\rho ^{t } \frac{W^{ (\gamma) \prime}(X_s-I_{s})-\Phi(\gamma)W^{ (\gamma)}(X_{s}-I_{s})}
{W^{(\gamma)}(X_{s }-I_{s})}  \,\sigma \, ds
\end{eqnarray*}
using  $B$. Second, consider the martingale
$$K_t: =\int_0^t \int_{(-1,0)} y \tilde{N}(dy,ds)$$
for the original filtration and let
\begin{eqnarray*}
\Gamma_t&=& -\int_0^{t\wedge \rho}   \frac{1}{1-Y_s} d\langle K,M\rangle_s + 1_{\{t\ge \rho\}} \int_\rho ^t \frac{1}{Y_s} d\langle K,M\rangle_s \\
& = & -\int_0^{t\wedge \rho} \frac{e^{-\Phi(\gamma)b}W^{(\gamma)}(b-I_{s})}{e^{-\Phi(\gamma)b}W^{(\gamma)}(b-I_{s})-e^{-\Phi(\gamma)X_{s}}W^{(\gamma)}
(X_{s}-I_{s})} \, d\langle K,M\rangle_s \\
& & + 1_{\{t\ge \rho\}} \int_\rho ^t  \frac{e^{-\Phi(\gamma)b}W^{(\gamma)}(b-I_{s})}{e^{-\Phi(\gamma)X_{s}}W^{(\gamma)}
(X_{s}-I_{s})}  d\langle K,M\rangle_s \\
&=&-\int_0^{t\wedge \rho} \int_{(-1,0)} y\, \frac{e^{-\Phi(\gamma)(X_{s}+y)}W^{ (\gamma)}(X_{s}+y-I_{s}) - e^{-\Phi(\gamma)X_{s}}W^{ (\gamma)}(X_{s}-I_{s})} {e^{-\Phi(\gamma)b}W^{(\gamma)}(b-I_{s})-e^{-\Phi(\gamma)X_{s}}W^{(\gamma)}
(X_{s}-I_{s})}   \,\Pi(dy) \, ds   \\
 & &
+ 1_{\{t\ge \rho\}} \int_\rho ^{t} \int_{(-1,0)}  y\,  \frac{e^{-\Phi(\gamma)y}W^{ (\gamma)}(X_{s}+y-I_{s}) - W^{ (\gamma)}(X_{s}-I_{s})} {W^{(\gamma)}(X_{s }-I_{s})}          \,\Pi(dy) \, ds
\end{eqnarray*}
for $t< \tau_b^+$. Finally, consider $N_t:=\int_0^t\int_{-\infty}^{-1}N(dy,ds)$ which is a Poisson process with mean $ct$ where $c:=\pi(-\infty,-1]$. The centered Poisson process $N_t-ct$  is  a martingale since $c$ is finite, that is, there are finitely many jumps with size less than or equal to $-1$. Then, we let
\begin{eqnarray*}
\Delta &=& -\int_0^{t\wedge \rho}   \frac{1}{1-Y_s} d\langle N,M\rangle_s + 1_{\{t\ge \rho\}} \int_\rho ^t \frac{1}{Y_s} d\langle N,M\rangle_s \\
 &=&-\int_0^{t\wedge \rho} \int_{-\infty}^{-1}   \frac{e^{-\Phi(\gamma)(X_{s}+y)}W^{ (\gamma)}(X_{s}+y-I_{s}) - e^{-\Phi(\gamma)X_{s}}W^{ (\gamma)}(X_{s}-I_{s})} {e^{-\Phi(\gamma)b}W^{(\gamma)}(b-I_{s})-e^{-\Phi(\gamma)X_{s}}W^{(\gamma)}
(X_{s}-I_{s})}   \,\Pi(dy) \, ds   \\
 & &
+ 1_{\{t\ge \rho\}} \int_\rho ^{t} \int_{-\infty}^{-1}    \frac{e^{-\Phi(\gamma)y}W^{ (\gamma)}(X_{s}+y-I_{s}) - W^{ (\gamma)}(X_{s}-I_{s})} {W^{(\gamma)}(X_{s }-I_{s})}          \,\Pi(dy) \, ds
\end{eqnarray*}

Now, let $$\bar{B}=B-\Lambda \; , \quad \quad \bar{K}=K-\Gamma \; , \quad \quad \bar{N}_t=N_t-ct -\Delta_t \; .$$
 It follows from \citep[Thm.VI.18]{prot} that $\bar{B},\bar{K}$ and $\bar{N}_t$  are martingales for the expanded filtration ${\cal F}^\rho$. Note that $\bar{N}$ is a counting process with compensator $ct+\Delta_t$ since $N_t -ct -\Delta_t$ is a martingale for ${\cal F}^\rho$, and define
\[
Q_t=\bar{B}_{\rho+t} - \bar{B}_\rho\; , \quad \quad L_t=\bar{K}_{\rho+t} - \bar{K}_\rho \;.
\]

Now, consider the post-$\rho$ process $X_{\rho+t} -X_\rho$ given by
\begin{equation} \label{postrho}
X_{\rho+t} -X_\rho=\mu \, t + \sigma(B_{\rho+t} - B_\rho)+ \int_\rho ^{\rho+t} \int_{(-1,0)} y \tilde{N}(dy,ds)+\int_\rho^{\rho+t} \int_{-\infty}^{-1} y N(dy,ds)
\end{equation}
Also, for $\rho+t < \tau^+_b $
\begin{eqnarray*}
\Lambda_{\rho+t}-\Lambda_\rho&=&  \int_\rho ^{\rho+t} \frac{W^{ (\gamma) \prime}(X_s-I_{s})-\Phi(\gamma)W^{ (\gamma)}(X_{s}-I_{s})}
{W^{(\gamma)}(X_{s }-I_{s})} \, \sigma \, ds
\end{eqnarray*}
and
\begin{eqnarray*}
\Gamma_{t+\rho}-\Gamma_\rho&=&
 \int_\rho ^{\rho+t} \int_{(-1,0)} y\,  \frac{e^{-\Phi(\gamma)y}W^{ (\gamma)}(X_{s}+y-I_{s}) - W^{ (\gamma)}(X_{s}-I_{s})} {W^{(\gamma)}(X_{s }-I_{s})}           \,\Pi(dy) \, ds  \\
\end{eqnarray*}
Putting $X_\rho= I_\rho$,  we can write (\ref{postrho}) as
\[
X_{\rho+t} -I_\rho= \mu  t + \sigma(\bar{B}_{\rho+t} - \bar{B}_\rho)+ \sigma (\Lambda_{\rho+t}-\Lambda_\rho) + \bar{K}_{\rho+t} - \bar{K}_\rho+\Gamma_{\rho+t}-\Gamma_\rho + \int_\rho^{\rho+t} \int_{-\infty}^{-1} y {N}(dy,ds)
\]
We see that $Q_t=\bar{B}_{\rho+t} - \bar{B}_\rho$ is a Brownian motion by L\'{e}vy characterization theorem,
and $L_t=\bar{K}_{\rho+t}-\bar{K}_\rho$ is a martingale with jumps independent from ${\cal F}^\rho$.  Therefore, the post-process $R_t:= X_{\rho+t}-X_\rho= X_{\rho+t}-I_\rho$, $0 < t < \tau_b^+ -\rho$,  satisfies the stochastic integral equation
\begin{eqnarray}
R_t &=&  \mu  t +\sigma Q_t+L_t
+ \, \sigma\int_0^t \frac{W^{ (\gamma) \prime}(R_s)-\Phi(\gamma)W^{ (\gamma)}(R_s)}
{W^{(\gamma)}(R_s)}  \, ds  +\int_0^t \int_{-\infty}^{-1} y \bar{N}(dy,ds) \nonumber \\
&&+\int_0^{t} \int_{(-1,0)}y \frac{e^{-\Phi(\gamma)y}W^{ (\gamma)}(R_s+y) - W^{ (\gamma)}(R_s)} {W^{(\gamma)}(R_s) }      \,\Pi(dy) \, ds  \label{postinf}
\end{eqnarray}
due to the fact that $X_s-I_s = X_s - I_\rho=R_{s-\rho}$ for $\rho < s < \tau_b^+ $, where $\bar{N}(dy,ds)$ now refers to $N(dy, ds+\rho)$, which forms a point process independent from ${\cal F}^\rho_\rho$ with compensator characterized with $\Delta$ above.
It can be shown that  \eqref{postinf} has a strong solution under the same conditions that the L\'{e}vy process $X$ exists, similar to the arguments given in the proof of Lemma \ref{Y} \citep{appl}. We will characterize the distribution of $R$ with its  infinitesimal generator. First, we have by Ito's formula  \citep[pg.66]{iked}
\begin{eqnarray*}
d\,F(R_t)&=&\left[\mu +\sigma \frac{W^{ (\gamma) \prime}(R_t)-\Phi(\gamma)W^{ (\gamma)}(R_t)}
{W^{(\gamma)}(R_t)}\right] F'(R_t) \, dt + \frac{\sigma^2}{2} F''(R_t)\, dt\\
&&    +\int_{(-1,0)}y \frac{e^{-\Phi(\gamma)y}W^{ (\gamma)}(R_t+y) - W^{ (\gamma)}(R_t)} {W^{(\gamma)}(R_t) }      \,\Pi(dy) \, dt                            \\
&& +\sigma F'(R_t) \, d Q_t + \int_{-\infty}^{-1} [F(R_{t^-}+y)-F(R_t)] \,  \bar{N}(dy, dt)  \\ && + \int_{(-1,0)} [F(R_{t^-}+y)-F(R_t)]  \, \bar{L}(dy, dt) \\
&& + \int_{(-1,0)} [F(R_t+y)-F(R_t)-F'(R_t)\, y ] \, \Pi(dy) \, dt    \\
&& \hspace{-2cm}+ \int_{(-1,0)} [F(R_t+y)-F(R_t)-F'(R_t)\, y ] \, \frac{e^{-\Phi(\gamma)y}W^{ (\gamma)}(R_t+y) - W^{ (\gamma)}(R_t)} {W^{(\gamma)}(R_t) }      \,\Pi(dy) \, dt
\end{eqnarray*}
for   $F\in C^2_b$, where we introduced the notation $\bar{L}(dy,ds)$ for $y\tilde{N}(dy, ds+\rho)$ where $\tilde{N}$ has the compensator characterized with $\Gamma$ with respect to ${\cal F}^\rho$.
We can now write the infinitesimal generator as
\begin{eqnarray*}
{\mathcal L}F(x) & =&  \left[\mu +\sigma \frac{W^{ (\gamma) \prime}(x)-\Phi(\gamma)W^{ (\gamma)}(x)}
{W^{(\gamma)}(x)}\right] F'(x)+\frac{\sigma^2}{2} F''(x) \\
& & + \int_{(-1,0)}y \frac{e^{-\Phi(\gamma)y}W^{ (\gamma)}(x+y) - W^{ (\gamma)}(x)} {W^{(\gamma)}(x) }      \,\Pi(dy) \\
&& + \int_{-\infty}^0 [F(x+y)-F(x)-F'(x)\, y\, 1_{\{y>-1\}} ] \, \Pi(dy)     \\
&& \hspace{-1cm}+ \int_{-\infty}^0 [F(x+y)-F(x)-F'(x)\, y\, 1_{\{y>-1\}} ] \, \frac{e^{-\Phi(\gamma)y}W^{ (\gamma)}(x+y) - W^{ (\gamma)}(x)} {W^{(\gamma)}(x) }      \,\Pi(dy)
\end{eqnarray*}
\finish

\label{sectionprev}

\section{Path decomposition conditioned on extremes}\label{sec5}

In this section, we display the distribution of various parts of the path conditioned on the supremum and the infimum based on the analysis of the previous sections.

\begin{thm}\label{thm2} Conditionally on $H_I<H_S$, $I_T=a$, $S_T=b$, it follows that
\begin{enumerate}
   \item the intermediate process $\{X_{{H_I}+t} :0\le t \le H_S -H_I \}$ is identical in law with $\{a+\bar{X}_{t} : 0\leq t\leq \tau^+_{b-a}\}$ where the infinitesimal generator of $\bar{X}$ is
   \begin{eqnarray*}
{\mathcal L}F(x) & =&  \left[\mu +\sigma \frac{\bar{W}^{ (\gamma) \prime}(x)-\Phi(\gamma)\bar{W}^{ (\gamma)}(x)}
{\bar{W}^{(\gamma)}(x)}\right] F'(x)+\frac{\sigma^2}{2} F''(x)  \\
& & + \int_{(-1,0)}y \frac{e^{-\Phi(\gamma)y}\bar{W}^{ (\gamma)}(x+y) - \bar{W}^{ (\gamma)}(x)} {\bar{W}^{(\gamma)}(x) }      \,\Pi(dy)  \\
&& + \int_{-\infty}^0 [F(x+y)-F(x)-F'(x)\, y\, 1_{\{y>-1\}} ] \, \Pi(dy)    \\
&& \hspace{-1cm}+ \int_{-\infty}^0 [F(x+y)-F(x)-F'(x)\, y\, 1_{\{y>-1\}} ] \, \frac{e^{-\Phi(\gamma)y}\bar{W}^{ (\gamma)}(x+y) - \bar{W}^{ (\gamma)}(x)} {\bar{W}^{(\gamma)}(x) }      \,\Pi(dy)
\end{eqnarray*} with the $\gamma-$scale function satisfying
\[\int\limits_0^\infty  {{e^{ - \lambda x}}} \bar{W}^{(\gamma)}(x)dx =\frac{1}{\bar{\psi}(\lambda)-\gamma}  =\frac{1}{\psi(\lambda + \Phi(\gamma))-2\gamma} ~~~~\mbox{for}~~~~\lambda\geq \Phi(2\gamma)-\Phi(\gamma)  \]
where $\Phi$ denotes the right inverse of $\psi.$

  \item the transition semigroup of the post-$H_S$ process is the $h$-transform of the law of a spectrally positive L\'{e}vy process killed at $T \wedge \hat{\tau}^+_{b-a} \wedge \hat{\tau}_0^-$,
      with
 \[
      h(z)=1-Z^{(\gamma)}(-z+b-a)+Z^{(\gamma)}(b-a)\frac{W^{ (\gamma)}(-z+b-a)}{W^{(\gamma)}(b-a)}-\frac{W^{ (\gamma)}(-z+b-a)}{W^{(\gamma)}(b-a)} \; ,
      \]
      that is,
      \[
P_t(x,dy)=\frac{h(b-y)}{h(b-x)}\, \p _{b-x}  \{  \hat{X}_t\in b-dy, t<T \wedge \hat{\tau}^+_{b-a} \wedge \hat{\tau}_0^- \}
\]
for $a<x,y<b$, and the entrance law is obtained as $x\to b-$.
\end{enumerate}
\end{thm}

\prf \textbf{i)}Conditioning on $S_{T}=b,$ the pre-supremum process is a spectrally negative L\'evy process with Laplace exponent
\[
\bar{\psi}(\lambda)=\psi(\lambda + \Phi(\gamma))-\gamma
\]
for $\lambda\ge -\Phi(\gamma)$, killed at the first passage time above $b,$ as given in Thm.\ref{thm3} and it is independent from the post-supremum process.
Comparing all the given conditions of Thm.\ref{thm2} with those in Prop.\ref{prop2}, we see that by definition of $\rho$, $X_{\rho+t}-X_{\rho}$ corresponds to the part between the last infimum, which is $a,$ and the level $b$ of the spectrally negative L\'evy process with the Laplace exponent $\bar{\psi}$. Therefore, the result follows from Prop.\ref{prop2}.\\
\textbf{ii)}
At an intermediate time $t<T$, the law of the post-$H_S$ process given $I_T=a$, $S_T=b$, $H_I<H_S$ is determined by transition semigroups $P_t(x,dy)$, $x<b$, $a<y<b$. This is equal to transition probability density at $b-y$ of a spectrally positive L\'evy process $\hat{X}$ killed at an exponential time $T$ and conditioned to stay positive and not go above $b-a$. Similar to post-$H_I$ process, the probability law of this process is written as
\begin{eqnarray}
& & \hspace{-1cm}\p_{b-x}\{\hat{X}_t\in b-dy, t<T\,|\, T< \hat{\tau}^+_{b-a} \wedge \hat{\tau}_0^-\} \nonumber \\
& = & \p_{b-x}\{\hat{X}_t\in b-dy, t<T, t<\hat{\tau}^+_{b-a} \wedge\hat{\tau}_0^-| \,T< \hat{\tau}^+_{b-a} \wedge \hat{\tau}_0^-\} \nonumber \\
&=& \frac{\p_{b-x}\{\hat{X}_t\in b-dy, t<T\wedge\hat{\tau}^+_{b-a} \wedge\hat{\tau}_0^-\}}{\p _{b-x}(T< \hat{\tau}^+_{b-a} \wedge \hat{\tau}_0^-)}\nonumber \\
& = & \frac{  \p_{b-y} \{ T<\hat{\tau}^+_{b-a} \wedge\hat{\tau}_0^-\} }{\p _{b-x}\{T< \hat{\tau}^+_{b-a} \wedge \hat{\tau}_0^-\}} \p_{b-x}  \{ \hat{X}_t\in b-dy, t<T \wedge \hat{\tau}^+_{b-a} \wedge \hat{\tau}_0^- \} \nonumber
\end{eqnarray}
by Markov property and the fact that $T$ is memoryless, and $\hat{\tau}_a^-:=\{t\geq 0: \hat{X}_t\leq a\}$. It follows that the law of the post-$H_S$ process is the $h$-transform of a spectrally positive L\'evy process killed at time $T \wedge \hat{\tau}^+_{b-a} \wedge \hat{\tau}_0^-$ with
\begin{eqnarray}\label{h} h(z) &=& \p _z\{T< \hat{\tau}^+_{b-a} \wedge \hat{\tau}_{0}^-\}=\p _{-z}\{T< \tau^-_{a-b} \wedge \tau_{0}^+\}\\\nonumber
&=& \p _{-z}\{T< \tau^-_{a-b}, \tau^-_{a-b}<\tau_{0}^+\}+\p _{-z}\{T<\tau_{0}^+, \tau_{0}^+ < \tau^-_{a-b}\} \\\nonumber
&=& \p _{-z}\{\tau^-_{a-b}<\tau_0^+\}-\E _{-z}[e^{-\gamma\tau^-_{a-b}}1_{\{ \tau^-_{a-b}<\tau_0^+\}}] \\\nonumber
&&\;+ \p _{-z}\{\tau_0^+ < \tau^-_{a-b}\}-\E _{-z}[e^{-\gamma\tau_0^+} 1_{\{\tau_0^+ <\tau^-_{a-b}\}}]\\\nonumber
&=& 1-Z^{(\gamma)}(-z+b-a)\\ \nonumber
&&\;+ Z^{(\gamma)}(b-a)\frac{W^{ (\gamma)}(-z+b-a)}{W^{(\gamma)}(b-a)}-\frac{W^{ (\gamma)}(-z+b-a)}{W^{(\gamma)}(b-a)} \end{eqnarray}
The $h$-function in \eqref{h} has been obtained using the expressions of the two sided exit formulas given in \citep[pgs. 17, 24]{KuztKypr}.
\finish

Recall that the maximum loss at time $t>0$ is defined by
\[
M_t^ - : = \mathop {\sup }\limits_{0 \le u \le v \le t} (X_u^{} - X_v^{}) \; ,
\]
\begin{cor} \label{cor2}
The distribution of the maximum loss of the post-supremum process, $M_{H_{S},T}^-$, when $H_I<H_S$, $I_T=a$, $S_T=b$ is given by
\begin{eqnarray*}&&\p(M_{H_{S},T}^-<d|H_I<H_S,I_T=a,S_T=b)\\
&=& 1- \left[1-Z^{(\gamma)}(b-a-d)+(Z^{(\gamma)}(b-a)-1)\frac{W^{ (\gamma)}(b-a-d)}{W^{(\gamma)}(b-a)}\right] \\
&& \;\;\;\cdot\;\frac{\gamma W^{(\gamma)}(d)-Z^{(\gamma)}(d)\frac{W^{ (\gamma)'}(d)}{W^{(\gamma)}(d)}}{-\gamma W^{(\gamma)}(b-a)+(Z^{(\gamma)}(b-a)-1)\frac{W^{(\gamma)'}(b-a)}{W^{(\gamma)}(b-a)}}
\end{eqnarray*}
for $0<d<b-a.$
\end{cor}
\prf
 We have
\begin{eqnarray*}\lefteqn {\p\{M_{H_{S},T}^-<d|H_I<H_S,I_T=a,S_T=b\}}\\
&= \displaystyle{ 1-\lim_{x \to b}\frac{h(d)}{h(b-x)}\, \p_{b-x}  \left\{ \hat{\tau_{d}}^+<T,\hat{\tau}_{d}^+< \hat{\tau}_{0}^-\right\}} \\
&= \displaystyle{1-\lim_{x \to b}\frac{h(d)}{h(b-x)}\, \p_{x-b}\left\{\tau_{-d}^-<T,\tau_{-d}^-<\tau_{0}^+\right\}}\\
&= \displaystyle{1-h(d)\lim_{x \to b}\frac{Z^{(\gamma)}(x-b+d)-Z^{(\gamma)}(d)\frac{W^{ (\gamma)}(x-b+d)}{W^{(\gamma)}(d)}}{1-Z^{(\gamma)}(x-a)+Z^{(\gamma)}(b-a)\frac{W^{ (\gamma)}(x-a)}{W^{(\gamma)}(b-a)}-\frac{W^{ (\gamma)}(x-a)}{W^{(\gamma)}(b-a)}}   }\\
&= \displaystyle{ 1- h(d)\frac{Z^{(\gamma)'}(d)-Z^{(\gamma)}(d)\frac{W^{ (\gamma)'}(d)}{W^{(\gamma)}(d)}}{-Z^{(\gamma)'}(b-a)+Z^{(\gamma)}(b-a)\frac{W^{ (\gamma)'}(b-a)}{W^{(\gamma)}(b-a)}-
\frac{W^{ (\gamma)'}(b-a)}{W^{(\gamma)}(b-a)}}  }
\end{eqnarray*}
which yields the result.
\finish

For removing the conditions in Thm.\ref{thm2} and Cor.\ref{cor2}, we provide the joint distribution of the supremum and the infimum. The  distribution is characterized in view of the formulas for the Laplace transform of the first passage times.
For $a<0<b,$ we have
\begin{eqnarray*} 
 \!\!  \p_{0}\{a<I_{T},S_{T}<b\} & = & \p_0\{T<\tau^{-}_a\wedge\tau^{+}_b\} \\
  &=& \p_{0}\{T<\tau^{-}_a,\tau^{-}_a<\tau^{+}_b\}+\p_{0}\{T<\tau^{+}_b,\tau^{+}_b<\tau^{-}_a\} \nonumber\\
      &=& 1-(\p_{0}\{T>\tau^{-}_a,\tau^{-}_a<\tau^{+}_b\}+\p_{0}\{T>\tau^{+}_b,\tau^{+}_b<\tau^{-}_a\})\nonumber\\
&=&1-(\E_{0}[e^{-\gamma\tau^{-}_a}1_{\{\tau^{-}_a<\tau^{+}_b\}}]+\E_{0}[e^{-\gamma\tau^{+}_b}1_{\{\tau^{+}_b<\tau^{-}_a\}}])\nonumber\\
& =& 1-(\E_{-a}[e^{-\gamma\tau^{-}_0}1_{\{\tau^{-}_0<\tau^{+}_{b-a}\}}]+\E_{-a}[e^{-\gamma\tau^{+}_{b-a}}1_{\{\tau^{+}_{b-a}<\tau^{-}_0\}}]).\nonumber
\end{eqnarray*}
The expectations in the above equation are known as two-sided exit from below and two sided-exit from above a level, respectively. For a spectrally negative L\'evy process, two sided exit problems are well defined and the expressions for them are well known \citep{KuztKypr}. After substituting these expressions, we find the joint distribution as
\begin{eqnarray*}
\lefteqn{ \p_{0}\{a<I_{T},S_{T}<b\}} \\ & \displaystyle{ \quad \quad \quad = 1-\left[Z^{(\gamma)}(-a)-Z^{(\gamma)}(b-a)\frac{W^{(\gamma)}
(-a)}{W^{(\gamma)}(b-a)}+\frac{W^{(\gamma)}(-a)}{W^{(\gamma)}(b-a)}\right]}\nonumber\\
&\!\!\!\!\!\! \!\!\!\!\!\! \!\! \displaystyle{=1-Z^{(\gamma)}(-a)+[Z^{(\gamma)}(b-a)-1]\frac{W^{(\gamma)}(-a)}{W^{(\gamma)}(b-a)} }. \nonumber
\end{eqnarray*}
This formula for the joint distribution of the extremes of a spectrally negative L\'{e}vy process generalizes that for a Brownian motion with drift \citep{ceren}.

Our motivation for path decomposition has been the ultimate aim of finding  the joint distribution of the maximum loss and the maximum gain, which result exists for standard Brownian motion as shown in  \citep{salm}.  The Brownian path is decomposed at the hitting times of the infimum and the supremum by considering  two cases, namely,
the infimum is attained before the supremum and vice versa.
In the first case, the maximum gain boils down to the range of the path, and in the second, the maximum loss is the range $S_T - I_T$.
Using the independence of the segments, the symmetry property of standard Brownian motion and the joint distribution of the infimum and the supremum, the joint distribution of the maximum gain and the maximum loss is found for both cases. Combining these yields the unconditional joint distribution.
For the spectrally negative L\'{e}vy process, the characterization of the pre-infimum part remains as an open problem. If solved, the result of this paper can be used towards the joint law of the maximum loss and the maximum gain as in \citep{salm} in view of duality between $X$ and $-X$.



\vspace{0.5cm}
\noindent{\bf Acknowledgement.}
We are grateful to Andreas E. Kyprianou for his
helpful discussions on the post-infimum process.
\vspace{0.5cm}
\\
\bibliographystyle{harvard}

\end{document}